# Why we shouldn't fault Lucas and Penrose for continuing to believe in the Gödelian argument against computationalism

**&**

**how we should, instead, use Gödel's reasoning to define *logical satisfaction*, *logical truth*, *logical soundness*, and *logical completeness* verifiably, and unarguably**

Bhupinder Singh Anand[1]

The only fault we can fairly lay at Lucas' and Penrose's doors, for continuing to believe in the essential soundness of the Gödelian argument, is their naïve faith in, first, non-verifiable assertions in standard expositions of classical theory, and, second, in Gödel's unvalidated interpretation of his own formal reasoning. We show why their faith is misplaced in both instances.

## 1. Introduction

Although most reasoned critiques (such as, for instance, [Bo90], [Br00], [Da93], [Fe96], [La98], [Le69], [Le89], [Pu95]) of Lucas' and Penrose's arguments against computationalism are unassailable, they do not satisfactorily explain why Lucas and Penrose - reasonable men, both - remain convinced of the essential soundness of their own arguments ([Lu96], [Pe96]).

[1] The author is an independent scholar. E-mail: re@alixcomsi.com; anandb@vsnl.com. Postal address: 32, Agarwal House, D Road, Churchgate, Mumbai - 400 020, INDIA. Tel: +91 (22) 2281 3353. Fax: +91 (22) 2209 5091.

A less technically critical review of their arguments is, indeed, necessary to appreciate the reasonability of their belief. It stems from the fact that, on the one hand, Lucas and Penrose have, unquestioningly, put faith in, and followed, standard expositions of classical theory in overlooking what Gödel has actually proven in Theorem VI ([Go31], p24) of his seminal 1931 paper [Go31] on formally undecidable arithmetical propositions; on the other, they have, similarly, put faith in, and uncritically accepted as definitive, Gödel's own, informal and unvalidated, interpretation of the implications of this Theorem ([Go31], p27).

They should not be taken to task on either count for their faith; it is standard expositions of Gödel’s reasoning that remain ambiguously silent on both issues.

In this paper, we show that such ambiguity and silence has been, and continues to be, both, misleading and unnecessary. We show, specifically, that in [Go31], Gödel has - albeit implicitly[2], and perhaps unwittingly and unknowingly - defined the *logical satisfaction* and *logical truth* of the formulas of an Arithmetic, the *logical soundness* of the Arithmetic itself, and the *logical completeness* of the Arithmetic, in an effective and verifiable manner within the Arithmetic.

Why the definitions have not been explicitly recognised by, both, Gödel (at least in [Go65]), and standard expositions of his reasoning (such as, for instance, [Be59], [Bo03], [Ch98], [Da93], [Fe96], [Me64], [R087], [Ro36], [Ro39], [Sh67], [Sm92], [Wa64]) is a mystery.

[2] “The method of proof which has just been explained can obviously be applied to every formal system which, first, possesses sufficient means of expression when interpreted according to its meaning to define the concepts (especially the concept “provable formula”) occurring in the above argument; and, secondly, in which every provable formula is true. In the precise execution of the proof, which now follows, we shall have the task (among others) of replacing the second of the assumptions just mentioned by a purely formal and much weaker assumption.” ([Go31], p9)

The question we need to ask is, however: How differently would we have viewed Gödel's reasoning, and its consequences, had Gödel defined *logical satisfiability*, *logical truth*, *logical soundness* and *logical completeness* explicitly as below.

## 2. Can verifiable *logical truth* be formalised in Peano Arithmetic?

Now, standard expositions of Tarski's Theorem [Ta36] - to the effect that the set of Gödel numbers of the formulas of any first-order Peano Arithmetic, which are *intuitively true* in the standard model[3] of the Arithmetic, is not arithmetical - appear to implicitly suggest that a verifiable '*logical truth*' of the formulas of standard Peano Arithmetic, under an interpretation, cannot be formalised in the Arithmetic.

Accepting this, seeming, implication unquestioningly, Lucas and Penrose use it explicitly as an arguable cornerstone of the Gödelian argument[4] [Lu61][Pe90][Pe94].

However, the crucial point provable by Gödel's reasoning in [Go31], but one whose significance has been overlooked, both, by him as well as by standard expositions of his reasoning, is that such a conclusion is not just arguable - it is false.

---

[3] We define the "standard interpretation" of first-order Peano Arithmetic, PA, as (cf. [Me64], p107):

"... the interpretation in which

(*a*) the set of non-negative integers is the domain,
(*b*) the integer 0 is the interpretation of the symbol 0,
(*c*) the successor operation (addition of 1) is the interpretation of the ' function (i.e., of $f_1^1$),
(*d*) ordinary addition and multiplication are the interpretations of + and .,
(*e*) the interpretation of the predicate letter = is the identity relation."

In other words, the interpreted, arithmetical, relation $R(x)$ is obtained from the formula $[R(x)]$ of PA by replacing every primitive, undefined, symbol of PA in the formula $[R(x)]$ by an intuitively interpreted mathematical symbol (i.e. a symbol that is a shorthand notation for some, semantically well-defined, concept of classical mathematics) as in (*a*)-(*e*).

[4] Since it is not germane to the purely logical issue sought to be raised here, we shall not concern ourselves with the specific content, arguments and objections to the Gödelian argument.

## 3. Tarski's definitions of *satisfiability* and *truth*

Now, the standard definitions of the *satisfiability*, and *truth*, of the formulas of a formal language, say L, under a well-defined interpretation[5], say M, are due to Tarski [Ta36].

Thus, a formula[6] [$R(x)$] of L is defined as *satisfied* under M if, and only if, it's corresponding interpretation, say $R^*(x)$, *holds*[7] in M for any assignment of a value $s$ that lies within the range of the variable $x$ in M.

The formula [$(\mathrm{A}x)R(x)$] of L is, then, defined as *true* under the interpretation M if, and only if, [$R(x)$] is *satisfied* under M.

---

[5] The word "interpretation" may be used both in its familiar, linguistic, sense, and in a mathematically precise sense; the appropriate meaning is usually obvious from the context. Mathematically, we follow the following definition ([Me64], §2, p49):

"An *interpretation* consists of a non-empty set D, called the *domain* of the interpretation, and an assignment to each predicate letter $A_j^n$ of an $n$-place relation in D, to each function letter $f_j^n$ of an $n$-place operation in D (i.e., a function from $D^n$ into D), and to each individual constant $a_i$ of some fixed element of D. Given such an interpretation, variables are thought of as ranging over the set D, and ~, =>, and quantifiers are given their usual meaning. (Remember that an $n$-place relation in D can be thought of as a subset of $D^n$, the set of all $n$-tuples of elements of D.)"

We note that the interpreted arithmetical relation, $R(x)$, in the standard model M of a Peano arithmetic, P, is obtained from the formula [$R(x)$] of the formal system P by replacing every primitive, undefined symbol of P in the formula [$R(x)$] by an interpreted mathematical symbol (i.e. a symbol that is a shorthand notation for some, semantically well-defined, concept of classical mathematics). So the P-formula [$(\mathrm{A}x)R(x)$] interprets as the sentence $(\mathrm{A}x)R(x)$, and the P-formula [$\sim(\mathrm{A}x)R(x)$] as the sentence $\sim(\mathrm{A}x)R(x)$.

We also note that the meta-assertions "[$(\mathrm{A}x)R(x)$] is a true sentence under the interpretation M of P", and "$(\mathrm{A}x)R(x)$ is a true sentence of the interpretation M of P", are equivalent to the meta-assertion "$R(x)$ is satisfied for any given value of $x$ in the domain of the interpretation M of P" ([Me64], p51).

[6] We use square brackets to indicate that the expression within the brackets is to be treated as a syntactic string of formal symbols only, devoid of any semantic content.

[7] Tarski's definitions are mathematically significant only if we assume that, given any $s$ in M, we can effectively determine that $R^*(s)$ *holds instantiationally* in M. Where M is an intuitive interpretation, such determination is sought to be postulated by the Church and Turing Theses, albeit, *algorithmically*. In Appendix 1 we argue that there is no justification for such strong postulation, and that it needs to be weakened to *instantiational* determination only - the minimum requirement of Tarski's definitions.

Moreover, the formula [~(A$x$)$R$($x$)] of L is, further, defined as *true* under the interpretation M if, and only if, [(A$x$)$R$($x$)] is *not true* under M.

Clearly, mathematical *satisfaction* and *truth* is, thus, defined relative only to an interpretation.

The Gödelian argument, quite reasonably, therefore, attempts to draw philosophical conclusions from the meta-logical status of the *intuitive satisfaction* and *intuitive truth* given in mathematical reasoning to the formulas of Peano Arithmetic under its standard, *intuitive*, interpretation.

## 4. Definition of verifiable *logical satisfaction* and *logical truth*

However, if we take M to be an interpretation of L in L itself, then we have the, classically overlooked, formalisation of the concepts of verifiable, and unarguable, *logical satisfaction*, and *logical truth*, of the formulas [$R$($x$)] and [(A$x$)$R$($x$)] of L, respectively, in L, as:

> The formula [$R$($x$)] of L is defined as *logically satisfied* under L if, and only if [$R$($s$)] is *provable*[8] in L for any term [$s$] that can be substituted for the variable [$x$] in [$R$($x$)].
>
> The formula [(A$x$)$R$($x$)] of L is *logically true* in L if, and only if, [$R$($s$)] is *logically satisfied* in L.

[8] Gödel has shown in [Go31] that the concept of 'provability' can be defined constructively, and in a verifiable manner, in recursively defined languages by means of primitive recursive functions and relations.

## 5. Definition of verifiable *logical soundness*

If we, further, define *logical soundness* as the property that the axioms of a theory are satisfied in the theory itself, and that the rules of inference preserve *logical truth*, then, it follows that the theorems of any *logically sound* theory are *logically true* in the theory.

It is straightforward to verify that first-order Peano Arithmetic is, indeed, *logically sound*.

## 6. Gödelian propositions

Now, even if the formula $[(Ax)R(x)]$ is not provable in L, it would be *logically true* in L if, and only if, the formula $[R(s)]$ is always provable in L for every well-defined term $[s]$ of L that can be substituted for $[x]$ in $[R(x)]$.

An instance of such a Gödelian proposition is, precisely, what Gödel has proven in his Theorem VI ([Go31], p24) for Peano Arithmetic[9] - by constructing a formula, $[(Ax)R(x)]$, of PA that is, itself, unprovable in PA, even though, for any given numeral $[n]$, $[R(n)]$ is provable in PA.

So, Gödel has actually constructed a formally *unprovable* Arithmetical formula that is not only *intuitively true* in the standard, *intuitive*, interpretation of the Arithmetic, but which is also *logically true* in the Arithmetic in a verifiable, and intuitionistically unobjectionable, manner that leaves no room for dispute as to its 'truth' status.

Moreover, since the Arithmetic can be shown to be *logically sound* - again in a verifiable, and intuitionistically unobjectionable, manner - standard expositions of Gödel's meta-reasoning no longer need appeal to the (arguable) assumption that the Arithmetic is *intuitively sound* under the standard interpretation[10].

[9] Although Gödel's arguments were developed in an explicitly defined formal system P of Arithmetic based on Dedekind's formulation of the Peano Axioms, they hold also for standard first-order Peano Arithmetic.

Prima facie, removing the debatable elements of *intuitive soundness* and *intuitive truth* from the Gödelian argument - which is built around the above construction - and replacing them with verifiable definitions of *logical soundness*, *logical satisfaction*, and *logical truth*, should help place the argument against computationalism in better perspective.

## 7. Standard expositions should appeal to verifiable *logical satisfaction* and *logical truth*, not to *intuitive truth*

However, standard expositions of Gödel's reasoning - including Gödel's own - continue to overlook such verifiable, and intuitionistically unobjectionable, formalisations of the concepts of the *logical satisfaction*, and the *logical truth*, of the formulas of Peano Arithmetic in the Arithmetic itself.

Instead, they admit into the foundations of first-order Peano Arithmetic concepts of unverifiable (hence arguable), *intuitive satisfaction* and *intuitive truth* that are only Platonically conceivable in individual, *intuitive*, 'standard' interpretations of the Arithmetic.

Surely Lucas and Penrose should not be faulted for treating this as a necessary, and definitive, omission, and for continuing to believe in the essential soundness of the Platonic, individual, *intuitive* interpretations and consequences of their own reasoning!

---

[10] "The system F is said to be *sound* for a class S of sentences if whenever F proves phi with phi in S then phi is *true* in the structure N of natural numbers" [Fe96].

## 8. There are *no* non-trivial non-standard models of first order Peano Arithmetic

The significance of defining *logical satisfaction*, *logical truth*, *logical soundness*, and *logical completeness*, verifiably, may be far-reaching (Appendix B).

For instance, standard expositions of Gödel's reasoning seem to derive legitimacy for the existence of non-trivial[11], non-standard, models of Peano Arithmetic only from Gödel's unvalidated assertion[12] - at the end of his proof of Theorem VI ([Go31], p27) - which implicitly implies that, if [(A$x$)$R$($x$)] is a Gödelian proposition of PA, then the axiomatic addition of [~(A$x$)$R$($x$)] to PA will, first, not invite inconsistency, and, second, will yield a formal system, PA+[~(A$x$)$R$($x$)], with a non-trivial, non-standard, model, say M', of first-order Peano Arithmetic in which ~$R$($s$) holds for some $s$ in the domain of M' that is not a natural number.

However, this, again, is demonstrably false, if we require, reasonably, that consistency demand the extended theory remain *logically sound* in a verifiable manner.

For, by the above application of Tarski''s definitions to PA itself, it would falsely imply that:

> Since [~(A$x$)$R$($x$)] is provable in PA+[~(A$x$)$R$($x$)], it is *logically true* in PA+[~(A$x$)$R$($x$)]; hence [(A$x$)$R$($x$)] is *logically false* in PA+[~(A$x$)$R$($x$)], and so [$R$($n$)] is not provable for every numeral [$n$] in PA+[~(A$x$)$R$($x$)].

---

[11] We can define a trivial, non-standard, model of Peano Arithmetic as one obtained by adding to it, for example, a new individual constant [$b$], and corresponding axioms such as [$b$] ≠ [0], [$b$] ≠ [1], [$b$] ≠ [2], …, [$b$] ≠ [$n$], …. (cf. [Me64], p117).

[12] "If one adjoins Neg(17Gen$r$) to $\kappa$, then one obtains a consistent, but not an ω-consistent, class of FORMULAS $\kappa'$. $\kappa'$ is consistent, for, otherwise, 17Gen$r$ would be $\kappa$-provable. $\kappa'$ is however not ω-consistent, for, by virtue of ~$Bew_\kappa$(17Genr) and (15), we have ($x$)$Bew_\kappa$Sb(r 17|Z($x$)). On the other hand, of course, $Bew_{\kappa'}$ [Neg(17Gen$r$)] holds."

Clearly, PA+[~(A$x$)$R$($x$)] is not *logically sound*. We cannot, therefore, assume - as Gödel does - that it must be consistent if PA is consistent.

## 9. *Intuitive truth* admits arguable Platonic conclusions

So, perhaps, one should not be too harsh on the, mathematically questionable, philosophical conclusions that Penrose and Lucas draw from the assertion - implicitly endorsable by standard expositions of Gödel's reasoning - that [Lu96]:

> … in the case of First-order Peano Arithmetic there are Gödelian formulae (many, in fact infinitely many, one for each system of coding) which are not assigned truth-values by the rules of the system, and which could therefore be assigned either TRUE or FALSE, each such assignment yielding a logically possible, consistent system. These systems are random vaunts, all satisfying the core description of Peano Arithmetic.

## 10. Conclusions

What we have highlighted, above, is that the Gödelian argument draws sustenance from the fact that standard expositions of Gödel's reasoning appeal to the *intuitive satisfaction* and *intuitive truth* of the formulas of Peano Arithmetic, and the *intuitive soundness* of the Arithmetic, under an *intuitive* (hence arguable) standard interpretation; the concepts are not defined explicitly in an effectively verifiable manner.

We have shown that this ambiguous appeal, however, is easily avoided by using Gödel's reasoning to define the concepts of *logical satisfaction*, *logical truth*, and *logical soundness* in an effectively verifiable manner.

Moreover, if we define an Arithmetic as *logically complete* if, and only if, every *logically true* formula is provable in the Arithmetic, and require, further, that the *provable* *formulas* of a language must be *logically sound*, then we have the, more illuminating, interpretation of Gödel's reasoning as the assertion that Peano Arithmetic is not only *intuitively incomplete* - as asserted by, both, Gödel and standard expositions of his reasoning - but it is also *logically incomplete*.

Further, unlike the standard expositions of Gödelian incompleteness, which are rooted in the concept of an *intuitive*, hence arguable, *truth* in the standard model of Peano Arithmetic, *logical incompleteness* is not susceptible to the Gödelian argument.

## Appendix A: Commentary on mathematical objects and mathematical truth

The underlying issue, here, seems to be whether we can arrive at a common consensus with Lucas and Penrose on how we are to treat the terms 'mathematics' and 'computation'.

If we can begin by agreeing that, by 'mathematics', we mean those of our abstract mental concepts that can be expressed precisely and unambiguously in recursively definable mathematical languages, then the remaining issue is simply that of finding a mutually acceptable definition of 'computation' .

Now, there seems to be a curious, common, reluctance to highlight the fact that there are well-defined mathematical functions and relations that, classically, have been accepted as effectively computable / decidable, yet which are treated as 'uncomputable / undecidable' in current expositions of classical theory!

For instance, Dedekind's definition of a real number in terms of cuts, and the equivalent definition in terms of Cauchy sequences are, both, effectively computable / decidable.

This is why Chaitin [Ch98] can claim that his Omegas define real numbers (assuming that the definitions are valid), since any digit of a given definition of an Omega can be effectively computed mechanically. However, there is no single algorithm that can effectively decide the value of any digit of a given Omega.

### A.1 Well-defined functions *can* be effectively computable *instantiationally*, but not *algorithmically*

The straightforward way of expressing this phenomenon should be to say that there are well-defined real numbers that are *instantiationally* computable, but not *algorithmically* computable.

So why is this terminology uncomfortable for current expositions of classical theory, and why should the Omegas, amongst other, similarly definable, functions be termed as 'uncomputable' even in Computability Theory[13]?

### A.2 We use the term 'exists' ambiguously

The deeper issue here seems to be that, when using language to express the abstract objects (elements) of our individual, and common, mental 'concept spaces', we use the word 'exists' loosely in three senses, without making explicit distinctions between them.

First, we may mean that an individually conceivable object exists within a language L if it lies within the range of the variables of L. The existence of such objects is necessarily derived from the grammar and rules of construction of the appropriate constant terms of the language - generally finitary in recursively defined languages - and can be termed as constructive in L by definition.

[13] "After all, although no Turing machine computes the function $d$, we were able to compute at least its first few values, For since, as we have noted, $f_1 = f_1 = f_1$ = the empty function we have $d(1) = d(2) = d(3) = 1$. And it may seem that we can actually compute $d(n)$ for any positive integer $n$ - if we don't run out of time." ([Bo02], Ch. 4, Uncomputability, p37)

Second, we may mean that an individually conceivable object exists under a formal interpretation of L in another formal language, say L', if it lies within the range of a variable of L under the interpretation.

Again, the existence of such an object in L' is necessarily derivable from the grammar and rules of construction of the appropriate constant terms of L', and can be termed as constructive in L' by definition.

Third, we may mean that an individually conceivable, object exists, in an interpretation M of L, if it lies within the range of an interpreted variable of L, where M is a Platonic interpretation of L in an individual's subjective mental conception (a la Brouwer).

Clearly, the debatable issue is the third case.

### A.3 Can we correlate diverse, individually conceivable, interpretations unambiguously?

So the question is whether we can - and, if so, how we may - correspond the, Platonically conceivable, objects of various individual interpretations of L, say M, M', M'', ..., unambiguously to the mathematical objects that are definable as the constant terms of L.

If we can achieve this, we can, then, attempt to relate L to a common external world, and try to communicate effectively about our individual mental concepts of the world that we accept as lying, by consensus, in a common, Platonic, 'concept space'.

## A.4 The central role of the standard interpretation of first-order Peano Arithmetic

For mathematical languages, such an intuitionistically unobjectionable, common, ‘concept space’, is, implicitly, accepted as the set of individual, *intuitive*, Platonically conceivable, perceptions - M', M'', M''', ... - of the definition of the standard, *intuitive*, interpretation, say M, of Dedekind’s formulation of the Peano Axioms.

Reasonably, if we intend a language, or a set of languages, to be adequate, first, for the expression of the abstract concepts of an individual consciousness, and, second, for the unambiguous and effective communication of those of such concepts that we can accept as lying within our common concept space, then we need to give effective guidelines for determining the, Platonically conceivable, mathematical objects of an individual perception of M that we can agree upon, by common consensus, as corresponding to the constants (mathematical objects) definable within the language.

## A.5 The role of Church’s and Turing’s Theses in legitimising the standard interpretation

Now, in the case of mathematical languages in standard expositions of classical theory, this role is sought to be filled by the Church and Turing Theses. Their standard formulations postulate that every effectively computable number-theoretic function (or relation, treated as a Boolean function) of M is partial recursive / Turing-computable.

However, curiously, even Computability Theory is reluctant to note that these Theses do not succeed in their objective completely.

Thus, even if we accept the standard formulations of the Theses, we still cannot conclude that we have specified explicitly that the domain of M consists of only constructive mathematical objects that can be represented in the most basic of our formal mathematical languages, namely, first-order Peano Arithmetic and Recursive Arithmetic.

**A.6 The standard formulation of CT violates the principle of Occam's razor**

The reason seems to be that the Church and Turing Theses - CT for short - are postulated as strong identities[14], which, prima facie, go beyond the minimum requirements for the correspondence between the, Platonically conceivable, mathematical objects of M and those of PA and Recursive Arithmetic.

This violation of the principle of Occam's Razor is highlighted if we note that every recursive function (or relation) is not identical to a unique arithmetical function (or relation), but only *instantiationally* equivalent to an infinity of arithmetical functions (or relations)[15].

Thus, the standard form of CT only postulates as constructive the algorithmically computable number-theoretic functions of M.

It leaves open the question of the significance that we are to permit to the individual, Platonically conceivable, non-constructive, elements of M.

It also obscures the issue of whether there are constructive, *instantiationally* computable but *algorithmically* uncomputable, number-theoretic functions and relations.

---

[14] For instance, in [Ch36], Church writes: "We now define the notion, already discussed, of an effectively calculable function of positive integers by identifying it with the notion of a recursive function of positive integers."

Correspondingly, Turing notes, in [Tu36], that: "The theorem that all effectively calculable sequences are computable and its converse are proved below in outline".

[15] See Theorem VII in [Go31].

## A.7 Standard expositions of classical theory are restrictive because they identify *effective* computability with *algorithmic* computability

Thus, standard expositions of classical theory imply - albeit implicitly - that only *algorithmically* computable functions (and relations) can be termed as constructive.

It is the tacit acceptance of this implicit implication that prevents, for instance, a constructive definition of what we are individually able to conceive as a random real number - one whose digits are *instantiationally* computable effectively, but which are not computable *algorithmically* (such as, for instance, Chaitin's Omegas, assuming that they are, indeed, well-defined real numbers).

## A.8 We *can* define arithmetical truth effectively

Now, such implicit implication can be avoided by explicitly recognising that there are well-defined mathematically expressible functions (and relations) which can, intuitively, be termed as effectively computable / decidable *instantiationally* (i.e., in any given instance, by some mechanical method that depends on the given instance), even though they are not computable / decidable *algorithmically* (i.e., by a common mechanical method that applies to every given instance).

Recognition of this immediately allows us to define *logical truth* effectively (at least for Peano Arithmetic, as outlined earlier), under Tarski's definitions of the *satisfiability* and *truth* of the formal expressions of a language under a well-defined interpretation - despite the perceived limitations of Gödel's 'Incompleteness' Theorems, Tarski's Theorem that the set of Gödel numbers of *true* arithmetical statements is not arithmetical, Turing's Halting Theorem, and Cantor's diagonal argument.

## A.9 Gödel's 'Incompleteness' Theorems, Tarski's Theorem, Turing's Halting Theorem, and Cantor's diagonal argument in perspective.

For instance, the reason there is a formally *unprovable* expression of Peano Arithmetic that, under the standard interpretation of the language, translates, both, as *intuitively*, and as *logically*, true - under Tarski's definitions of the *satisfiability* and *truth* of such expressions under an interpretation - simply means that even standard formulations of Tarski's definitions admit the possibility of arithmetical relations as *true* that are *instantiationally*, but not *algorithmically*, *true* for any given set of natural number values in the interpretation, whereas a *provable* expression of PA necessarily translates as a relation that is *algorithmically* decidable as *true* in the interpretation.

So, Gödel's *unprovable*, but *intuitively true*, arithmetical relation may simply be one that is effectively decidable as *logically true* in every instance (which Gödel has proved), but which may not be *algorithmically* decidable as *logically true* (a point that does not seem to have been considered explicitly in current expositions of classical theory).

This, in essence, would be a reasonable interpretation of Tarski's Theorem - that there are arithmetical relations that are effectively decidable *instantiationally*, but not *algorithmically*.

Further, since we can define an arithmetical relation as effectively decidable *instantiationally* if, and only if, each instance of it were *provable* in Peano Arithmetic, it would be more illuminating to say that Peano Arithmetic can be termed as *instantiationally* complete, but not *algorithmically* complete.

Similarly, Turing's Halting Theorem (as also Cantor's diagonal construction) can be interpreted as asserting that there are well-defined number-theoretic functions (real numbers) that are effectively computable *instantiationally*, but not *algorithmically*.

**A.10 Consequences of strengthening the weakened forms of CT with a plausible arithmetical Provability Thesis**

Moreover, under an intuitionistically unobjectionable - and plausible - Provability Thesis (to the effect that a total arithmetical relation is PA-*provable* if, and only if, it is *algorithmically* decidable), we can define an architecture for a trio of Turing machines such that it can define a total number-theoretic function that is effectively computable *instantiationally*, but not *algorithmically*[16].

Obviously, acceptance of such a Thesis as intuitionistically unobjectionable implies that the standard forms of the Church and Turing Theses need to be weakened to equivalences.

It also implies that there are no non-trivial, non-standard, models of Peano Arithmetic, an immediate corollary of which is that P ≠ NP (Appendix B).

However, a more interesting benefit of weakening the Church and Turing Theses is that they, then, provide the equivalent *instantiational*, arithmetical, *completeness* - in first-order theory - that is provided in second-order Peano Arithmetic (which formalises our concepts of the natural numbers as expressed by Dedekind's Peano Postulates) by the second-order Induction Axiom.

Such *completeness* would be of significance to the computational thesis - which Lucas and Penrose attempt to refute by the Gödelian argument - that all 'mathematics' is precisely captured by 'computation'.

[16] The author develops this argument further in various arXived, but unpublished, papers.

## Appendix B: If first-order Peano Arithmetic has no consistent non-trivial, non-standard, models, then P ≠ NP

P ≠ NP is the central open problem in complexity theory [Co00], one of whose formulations is the following [Ra02]:

> "Is there a polynomial time algorithm A that gets as input a Boolean formula *f* and outputs 1 if and only if *f* is a tautology? P ≠ NP states that there is no such algorithm."

Now, Gödel has defined a formula, $[R(x)]$, such that:

(*i*) $[R(x)]$ is constructible in standard, first-order, Peano Arithmetic, PA;

(*ii*) we can prove, meta-mathematically, that $[R(x)]$ translates as an arithmetical tautology, $R(x)$, under the standard interpretation of the Arithmetic;

(*iii*) $[R(x)]$ is not provable in the Arithmetic.

The question arises: Is $R(x)$ Turing-decidable as TRUE?

If we assume, first, the thesis that every total arithmetical relation that is Turing-decidable as TRUE is PA-provable, then $R(x)$ is not Turing-decidable as TRUE, and, so, P ≠ NP.

If we assume, however, that there is a total arithmetical relation that is Turing-decidable as TRUE, but which is not PA-provable, then this implies that there is a consistent, non-trivial, non-standard, model of PA, in which $[R(s)]$ is satisfied for some term $[s]$ of the interpretation that is not a natural number.

We conclude that, if PA has no consistent, non-trivial, non-standard models, then, under the above expression of the PvNP problem, P ≠ NP.

## References

[Be59] Beth, E. W. 1959. *The Foundations of Mathematics*. Amsterdam. North-Holland

[Bo90] Boolos, G. S., et al. 1990. *An Open Peer Commentary on The Emperor's New Mind*. In Behavioral and Brain Sciences 13 (4) (1990) 655.

[Bo03] Boolos, G. S., Burgess, J. P., Jeffrey R. C. 2003. Computability and Logic (4th ed). Cambridge University Press, Cambridge.

[Br00] Bringsjord, S., Xiao, H. 2000. *A Refutation of Penrose's Gödelian Case Against Artificial Intelligence*. Journal of Experimental and Theoretical Artificial Intelligence 12: 307–329.

Offprint available at http://kryten.mm.rpi.edu/refute.penrose.pdf

[Ch36] Church, A. 1936. *An unsolvable problem of elementary number theory*. Am. J. Math., Vol. 58, pp. 345-363. Also in M. Davis (ed.). 1965. The Undecidable. Raven Press, New York.

[Ch98] Chaitin, G. J. 1998. The Limits of Mathematics. Springer-Verlag, Singapore.

[Co00] Cook, S. *The P versus NP Problem*. 2000. Official description provided for the Clay Mathematical Institute, Cambridge, Massachusetts.

<http://www.claymath.org/millennium/P_vs_NP/Official_Problem_Description.pdf>

[Da93] Davis, M. 1993. *How subtle is Gödel's theorem? More on Roger Penrose*. Behavioral and Brain Sciences, 16, 611-612.

[Fe96] Feferman, S. 1996. *Penrose's Gödelian argument*. PSYCHE 2 (1996) 21-32.

[Go31] Gödel, Kurt. 1931. *On formally undecidable propositions of Principia Mathematica and related systems I*. Translated by Elliott Mendelson. In M. Davis (ed.). 1965. The Undecidable. Raven Press, New York.

[Go65] Gödel, Kurt. 1965. *On Undecidable Propositions of Formal Mathematical Systems*. In M. Davis (ed.). 1965. The Undecidable. Raven Press, New York. (Reprinted from mimeographed notes prepared by S.C.Kleene and J.B.Rosser on lectures given by Professor Kurt Gödel at the Institute for Advanced Study, Princeton, during the spring of 1934; prepared and edited specially for the anthology by Professor Kurt Gödel.)

[La98] LaForte, G., Hayes, P. J., Ford, K. M. 1998. *Why Gödel's theorem cannot refute computationalism*. Artificial Intelligence, 104:265-286, 1998.

[Le69] Lewis, D. 1969. *Lucas against mechanism*. Philosophy 44 (1969) 231-233.

[Le89] Lewis, D. 1989. *Lucas against mechanism II*. The Canadian J. of Philosophy 9 (1989) 373-376.

[Lu61] Lucas, John, R. 1961. *Minds, Machines and Gödel*. First published in *Philosophy*, XXXVI, 1961, pp.(112)-(127); reprinted in *The Modeling of Mind*, Kenneth M. Sayre and Frederick J. Crosson, eds., Notre Dame Press, 1963, pp.[269]-[270]; and *Minds and Machines*, ed. Alan Ross Anderson, Prentice-Hall, 1954, pp.{43}-{59}.

<*Web page*: http://users.ox.ac.uk/~jrlucas/Godel/mmg.html>

[Lu96] Lucas, John, R. 1996. *The Gödelian Argument: Turn Over the Page*. A talk given on 25/5/96 at a BSPS conference in Oxford.

<*Web page*: http://users.ox.ac.uk/~jrlucas/Godel/turn.html>

[Me64] Mendelson, Elliott. 1964. Introduction to Mathematical Logic. Van Norstrand, Princeton.

[Pe90] Penrose, R. (1990, Vintage edition). The Emperor's New Mind: Concerning Computers, Minds and the Laws of Physics. Oxford University Press.

[Pe94] Penrose, R. (1994). Shadows of the Mind: A Search for the Missing Science of Consciousness. Oxford University Press.

[Pe96] Penrose, R. 1996. *Beyond the Doubting of a Shadow. A Reply to Commentaries on Shadows of the Mind*. PSYCHE, 2(23), January 1996

<*Web page*: http://psyche.cs.monash.edu.au/v2/psyche-2-23-penrose.html>

[Pu95] Putnam, H. 1995. *Review of Shadows of the Mind*. In Bulletin of the American Mathematical Society 32 (1995) 370-373.

[Ra02] Raz, Ran. 2002. P ≠ NP, Propositional Proof Complexity, and Resolution Lower Bounds for the Weak Pigeonhole Principle. ICM 2002, Vol. III, 1–3.

[Ro87] Rogers, H. Jr. 1987. Theory of Recursive Functions and Effective Computability. The MIT Press, Cambridge, Massachsetts.

[Ro36] Rosser. J. B. 1936. *Extensions of some Theorems of Gödel and Church*. The Journal of Symbolic Logic. Vol. 1. (1936), pp.87-91.

[Ro39] Rosser. J. B. 1959. *An informal exposition of proofs of Gödel's Theorem and Church's Theorem*. The Journal of Symbolic Logic. Vol. 4. (1939), pp.53-60.

[Sh67] Shoenfield, J. R. 1967. Mathematical Logic. Published for The Association for Symbolic Logic by A. K. Peters, Ltd., Natick, Massachussetts.

[Sm92] Smullyan, R. 1992. Gödel's Incompleteness Theorems. Oxford University Press. New York.

[Ta36] Tarski, A. 1936. *Der Wahrheitsbegriff in den formalisierten Sprache*. Studia Philos., Vol. 1. Expanded English translation, *The concept of truth in the languages of the deductive sciences*, in *Logic, Semantics, Metamathematics, papers from 1923 to 1938* (p152-278), ed. John Corcoran. 1983. Hackett Publishing Company, Indianapolis.

[Tu36] Turing, Alan. 1936. *On computable numbers, with an application to the Entscheidungsproblem*. Proceedings of the London Mathematical Society, ser. 2. vol. 42 (1936-7), pp.230-265; corrections, Ibid, vol 43 (1937) pp. 544-546.

<*Web page*: http://www.abelard.org/turpap2/tp2-ie.asp - index>

[Wa64] Wang, Hao. 1964. *A Survey of Mathematical Logic*. Amsterdam : North Holland.

(*Created: Monday*, *12th June 2006*, *1:12:32 PM IST*, *by re@alixcomsi.com*)